\documentclass[12pt]{amsart}

\input xy
\xyoption{all}

\begin{document}

\title[Relative Completions and Cohomology]{Relative completions and the 
cohomology of linear groups over local rings}

\author{Kevin P. Knudson}\thanks{Partially supported by NSF grant no. 
DMS-0070119}
\address{Department of Mathematics, 1150 F/AB, Wayne State University, 
Detroit, Michigan 48202}
\email{knudson@math.wayne.edu}

\newtheorem{theorem}{Theorem}[section]
\newtheorem{prop}[theorem]{Proposition}
\newtheorem{lemma}[theorem]{Lemma}
\newtheorem{cor}[theorem]{Corollary}
\newtheorem{conj}[theorem]{Conjecture}

\newcommand{\zz}{{\mathbb Z}}
\newcommand{\zq}{{\mathbb Q}}
\newcommand{\zzp}{{\mathbb Z}_p}
\newcommand{\zzl}{{\mathbb Z}_\ell}
\newcommand{\U}{{\mathcal U}}
\newcommand{\K}{{\mathcal K}}
\newcommand{\G}{{\mathcal G}}
\newcommand{\T}{{\mathcal T}}
\newcommand{\cP}{{\mathcal P}}
\newcommand{\Ox}{{\mathcal O}_{X,x}}
\newcommand{\mx}{{\mathfrak m}_x}
\newcommand{\Ohat}{\hat{\mathcal O}_{X,x}}
\newcommand{\lslnk}{{\mathfrak sl}_n(k)}
\newcommand{\fp}{{\mathbb F}_p}
\newcommand{\fl}{{\mathbb F}_\ell}
\newcommand{\fpbar}{\overline{\mathbb F}_p}
\newcommand{\slnk}{SL_n(k)}
\newcommand{\slnkt}{SL_n(k[t])}
\newcommand{\slnpow}{SL_n(k[[T]])}
\newcommand{\slno}{SL_n({\mathcal O}_{X,x})}
\newcommand{\slnohat}{SL_n(\hat{{\mathcal O}}_{X,x})}
\newcommand{\slnz}{SL_n(\zz)}
\newcommand{\slnfp}{SL_n(\fp)}
\newcommand{\ra}{\rightarrow}
\newcommand{\lra}{\longrightarrow}
\newcommand{\bop}{\bigoplus}
\newcommand{\hcts}{H^2_{\text{cts}}}
\newcommand{\mhcts}{H^2_{\text{\em cts}}}

\maketitle

For a discrete group $G$ there are two well-known completions.  The 
first is the Malcev (or unipotent) completion.  This is a prounipotent group
$\U$, defined over $\zq$, together with a homomorphism $\psi:G\ra \U$
that is universal among maps from $G$ into prounipotent $\zq$-groups.
To construct $\U$, it suffices to consider the case where $G$ is nilpotent;
the general case is handled by taking the inverse limit of the Malcev
completions of the $G/\Gamma^rG$, where $\Gamma^\bullet G$ denotes
the lower central series of $G$.  If $G$ is abelian, then $\U=G\otimes\zq$.
We review this construction in Section \ref{unipotent}.

The second completion of $G$ is the $p$-completion.  For a prime $p$, we
set $G^{\wedge p} = \varprojlim G/\Gamma^r_p G$, where $\Gamma^\bullet_p G$ is the $p$-lower central series of $G$.  If $G$ is a
finitely generated abelian group, then $G^{\wedge p}=\zzp\otimes G$,
where $\zzp$ is the ring of $p$-adic integers \cite{bousfield}.  The
group $G^{\wedge p}$ is a pro-$p$-group and each $G/\Gamma^r_pG$ is
nilpotent provided $H_1(G,\fp)$ is finite dimensional.

Both of these completions are instances of a general construction.  Let $k$
be a field.  The unipotent $k$-completion of a group $G$ is a prounipotent
$k$-group $\U_k$ together with a homomorphism $G\ra\U_k$.  The group
$\U_k$ is required to satisfy the obvious universal mapping property.  The
Malcev completion is the case $k=\zq$ and the $p$-completion is the case
$k=\fp$.  This construction for other fields $k$ is probably well-known
to the experts, but it does not seem to be in the literature.

One reason to study such completions is that they may be used to gain 
cohomological information about the groups $G$ and $\U$.  Indeed, the 
restriction map 
$\hcts(\U,k)\ra H^2(G,k)$ is injective (the definition of $\hcts$ will be
recalled below).  This allows one to obtain either a lower bound for $\dim_k
H^2(G,k)$ or an upper bound for $\dim_k \hcts
(\U,k)$.

Unfortunately, the group $\U_k$ may be trivial (e.g., $G$ perfect, or
more generally, if $H_1(G,k)=0$).  To circumvent this, Deligne suggested
the notion of relative completion.  Suppose that $\rho:G\ra S$ is a
representation of $G$ in a reductive group $S$ defined over $k$.  Assume
that the image of $\rho$ is Zariski dense.  The completion of $G$ relative
to $\rho$ is a proalgebraic group $\G$ over $k$ which is an extension of 
$S$ by a prounipotent group $\U$ together with a map $G\ra\G$.  The group
$\G$ should satisfy the obvious universal mapping property.

The basic theory of relative completion in characteristic zero was worked
out by R.~Hain\cite{hain1}.  Many of his results remain valid in positive
characteristic.  We shall review this in Section \ref{relative} below.  For
examples of relative completions in characteristic zero beyond those
presented here, the reader is referred to Hain's study \cite{hain1} of the
completion of the mapping class group $\Gamma_{g,r}$ of a surface $S$
of genus $g$ with $r$ marked points relative to its symplectic representation as the group of automorphisms of $H_1(S,\zz)$.  Other 
examples, due to the author \cite{knudson2}, are the groups $SL_n(\zz[t])$
and $SL_n(\zz[t,t^{-1}])$ relative to their obvious representations in 
$\slnz$.  Recent work by Hain and M.~Matsumoto \cite{hainmat} tackles
the absolute Galois group $\text{Gal}(\overline{\zq}/\zq)$ relative to the
cyclotomic character $\chi:\text{Gal}(\overline{\zq}/\zq)\ra \zzl^\times$.

In this paper we study the completions of groups such as $\slno$, where
$\Ox$ is the local ring of a closed point $x$ on a smooth affine curve $X$.
We also study the completion of $\slnz$ relative to the reduction map
$\slnz\ra\slnfp$.  We use these completions to compute the second 
continuous cohomology of special linear groups over complete local rings.
A sample result is the following.

\medskip

\noindent {\bf Theorem \ref{powercoho}.} {\em If $k$ is a finite field
or a number field, then for all $n\ge 3$, $\mhcts(\slnpow,k)=0$.}

\medskip

As far as the author knows, this is the first calculation of continuous
cohomology with coefficients of the same characteristic as $k$.

This paper is organized as follows.  In Section \ref{construct}, we 
review Hain's construction of the completion of $\rho:G\ra S$.  In
Section \ref{unipotent}, we discuss in detail the case where $S$ is the
trivial group; this is the unipotent completion mentioned above.  In 
Section \ref{examples}, we present several examples of unipotent
completions.  Section \ref{relative} deals with the basic theory of 
relative completion and the computation of examples.  Finally, in Section
\ref{cohomology} we carry out some cohomology calculations.

\medskip

\noindent {\em Acknowledgements.}  I would like to thank Dick Hain
and Chuck Weibel for many useful discussions.  I am also grateful to
the referee for insisting that I not omit so many details.

\section{A Construction of the Relative Completion}\label{construct}

The following construction is due to R.~Hain \cite{hain1}.  Let $G$ be a
group and suppose that $\rho:G\ra S$ is a Zariski dense representation in
a reductive $k$-group $S$.  Consider all commutative diagrams of the 
form
$$\xymatrix{
1\ar[r] & U\ar[r] & E\ar[r] & S\ar[r] & 1 \\
            &            & G\ar[u]^{\tilde{\rho}}\ar[ur]_{\rho} &  & }$$
where $E$ is a linear algebraic group over $k$, $U$ is a unipotent subgroup
of $E$, and the image of $\tilde{\rho}$ is Zariski dense.  It is clear how to
define morphisms of such diagrams.  Moreover, the collection of such 
diagrams forms an inverse system \cite{hain1}, Prop. 2.1.  The completion
$\G$ of $G$ relative to $\rho$ is defined to be the inverse limit
$$\G=\varprojlim E$$ over all the above commutative diagrams.  The
kernel of the map $\G\ra S$ will be called the {\em prounipotent
radical} of $\G$.

The group $\G$ satisfies the following universal mapping property 
\cite{hain1}, Prop. 2.3.  Suppose that ${\mathcal E}$ is a proalgebraic
group over $k$ and that ${\mathcal E}\ra S$ is a homomorphism of
proalgebraic groups with prounipotent kernel.  If $\phi:G\ra {\mathcal E}$
is a map whose composition with ${\mathcal E}\ra S$ is $\rho$, then
there is a unique map $\Phi:\G\ra{\mathcal E}$ in the category of
proalgebraic $k$-groups such that the diagram
$$\xymatrix{
 {}  & \G\ar[dd]_{\Phi}\ar[dr] & {} \\
G\ar[ur]\ar[dr]_{\phi} &  &  S \\
 {}   & {\mathcal E}\ar[ur] & }$$ 
commutes.

\section{Unipotent Completion}\label{unipotent}
The unipotent completion of a group is the case where $S$ is the trivial
group.  The following construction of the unipotent completion is due to
D.~Quillen \cite{quillen}, although he considered only the characteristic
zero case.

Let $k$ be a field.  For a group $G$, denote by $J_k$ the augmentation
ideal of the group ring $kG$; i.e., $J_k$ is the kernel of the map 
$\varepsilon:kG\ra k$ defined by $\varepsilon(g)=1$.
The map $g\mapsto (g-1) + J^2_k$ induces an isomorphism
$$H_1(G,k)\stackrel{\cong}{\lra} J_k/J^2_k.$$
Denote by $kG^\wedge$ the completed algebra
$$kG^\wedge =\varprojlim kG/J^l_k;$$
it is a complete Hopf algebra ($kG$ is a Hopf algebra via $\Delta(g)=g
\otimes g$).  Denote the completion of $J_k$ by $\hat{J}_k$.

Consider the set of grouplike elements
$${\mathcal P}=\{x\in kG^\wedge : \Delta(x)=x\otimes x, \varepsilon(x)
=1\}.$$  This is a subgroup of the group of units of $kG^\wedge$.  Define
a filtration of ${\mathcal P}$ by
$$\cP^l = \cP\cap (1+\hat{J}^l_k)$$
and set $$\cP_l = \cP/\cP^{l+1}.$$ Observe that $[\cP^i,\cP^j]\subseteq
\cP^{i+j}$.

\begin{lemma}\label{plalg}
 If $H_1(G,k)$ is finite dimensional, then $\cP_l$ may be given the structure
of a linear algebraic group over $k$.
\end{lemma}

\begin{proof}  Since $H_1(G,k)$ is finite dimensional and the graded algebra
$$\text{Gr}^\bullet_{\hat{J}} kG^\wedge = \bop_{i\ge 0} \hat{J}_k^i/
\hat{J}_k^{i+1}$$  is generated by
$$\text{Gr}^1_{\hat{J}} kG^\wedge = \hat{J}_k/\hat{J}^2_k \cong H_1(G,k),$$
the quotient algebra $kG^\wedge/\hat{J}^{l+1}_k$ is a finite dimensional
$k$-algebra. Define a homomorphism $\varphi:\cP\ra \text{Aut}(kG^\wedge/
\hat{J}_k^{l+1})$ by $\varphi(x)(u)=xu$.  Since the elements of
$\cP^{l+1}$ act trivially on $kG^\wedge/\hat{J}_k^{l+1}$, we have an
induced embedding $\overline{\varphi}:\cP_l\ra \text{Aut}(kG^\wedge/
\hat{J}_k^{l+1})$.  This endows $\cP_l$ with the structure of an algebraic
group over $k$. Since $\cP_l$ consists of those $x$ satisfying
the polynomial conditions $\varepsilon(x)=1$ and $\Delta(x)=x\otimes x$,
we see that $\overline{\varphi}(\cP_l)$ is a closed subvariety of
$\text{Aut}(kG^\wedge/\hat{J}_k^{l+1})$.
\end{proof}

Recall that an algebraic group $U\subset GL(V)$ over $k$ is called
{\em unipotent} if there is a filtration 
$$V=V^0\supset V^1\supset\cdots \supset V^m\supset 0$$ of $V$ such
that each element of $U$ acts trivially on the quotients $V^i/V^{i+1}$.
We claim that each $\cP_l$ is unipotent. To see this, note that we
have a filtration of $kG^\wedge/\hat{J}_k^{l+1}$ by powers of $\hat{J}_k$:
$$kG^\wedge/\hat{J}_k^{l+1}\supset \hat{J}_k/\hat{J}_k^{l+1}
\supset\cdots \supset \hat{J}_k^l/\hat{J}_k^{l+1}\supset 0.$$
If $g\in\cP_l$, then $g-1\in\hat{J}_k$ and hence $(g-1)\hat{J}_k^i
\subseteq \hat{J}_k^{i+1}$.  Thus, if $x\in\hat{J}_k^i$, we have
$gx=x+v$ for some $v\in\hat{J}_k^{i+1}$; i.e., $g$ acts trivially
on
$$\hat{J}_k^i/\hat{J}_k^{i+1}\cong (\hat{J}_k^i/\hat{J}_k^{l+1})/ 
(\hat{J}_k^{i+1}/\hat{J}_k^{l+1}).$$
Thus, $\cP_l$ is a unipotent group over $k$.  As such, it is isomorphic
to a subgroup of some $U_n(k)$ (\cite{borel}, p.~87), the subgroup of
$GL_n(k)$ consisting of upper triangular matrices with 1's on the
diagonal.  In particular, each $\cP_l$ is nilpotent.

We shall need the following fact about nilpotent groups.

\begin{lemma}\label{nilpth1} Let $f:M\ra N$ be a homomorphism of nilpotent
groups such that the induced map $f_*:H_1(M,\zz)\ra H_1(N,\zz)$ is
surjective.  Then $f$ is surjective.
\end{lemma}

\begin{proof}  Denote by $\Gamma^\bullet M$ and $\Gamma^\bullet N$ the
lower central series of $M$ and $N$, respectively.  Consider the induced
map of associated graded algebras
$$\text{Gr}(f):\text{Gr}^\bullet_\Gamma M\lra \text{Gr}^\bullet_\Gamma N,$$
where $\text{Gr}^\bullet_\Gamma G = \bop_{i\ge 1} \Gamma^iG/\Gamma^{i+1}G$.
This algebra is a Lie algebra over $\zz$ with bracket induced by the
commutator in $G$.  Note that as an algebra $\text{Gr}^\bullet_\Gamma G$
is generated by $\text{Gr}^1_\Gamma G = G/\Gamma^2 G = H_1(G,\zz)$.
Consider the commutative diagram
$$\xymatrix{
(\text{Gr}^1_\Gamma M)^{\otimes l}\ar@{>>}[d]\ar@{>>}[r]^{f_*^{\otimes l}} &
(\text{Gr}^1_\Gamma N)^{\otimes l}\ar@{>>}[d] \\
\text{Gr}^l_\Gamma M\ar[r]^{\text{Gr}^l(f)} & \text{Gr}^l_\Gamma N.}$$
Since $f_*^{\otimes l}$ is surjective by hypothesis and the vertical
arrows are also surjective, we see that $\text{Gr}^l(f)$ is surjective.

Now, say that $\Gamma^{m+1} M =\{1\}$, but $\Gamma^m M\ne \{1\}$.
Then since $\text{Gr}^l(f)$ is surjective for all $l$, we must
have $\Gamma^{m+1}N/\Gamma^{m+2}N=0$; that is, $\Gamma^{m+1}N=\Gamma^{m+2}N$.
This can occur only if $\Gamma^{m+1}N=\{1\}$ (since $N$ is nilpotent).
Consider the commutative diagram
$$\xymatrix{
1\ar[r] & \Gamma^m M\ar[r]\ar@{>>}[d]^{\text{Gr}^m(f)} & \Gamma^{m-1}M
\ar[r]\ar[d] & \text{Gr}^{m-1} M\ar@{>>}[d]^{\text{Gr}^{m-1}(f)}\ar[r] & 0 \\
1\ar[r] & \Gamma^m N\ar[r] & \Gamma^{m-1}N\ar[r] & \text{Gr}^{m-1}N\ar[r] & 0
}$$
This shows that $\Gamma^{m-1}M$ surjects onto $\Gamma^{m-1}N$.
A simple induction then shows that $f:M\ra N$ is surjective.
\end{proof}

Now, let $G$ be a group with $H_1(G,k)$ finite dimensional.
We have an inclusion $G\ra kG$.  Composing with the completion map
$kG\ra kG^\wedge$ gives a map $G\ra kG^\wedge$.  The image of this
map lies in $\cP$ and hence for each $l$ there is a homomorphism
$G\ra \cP_l$.  We claim that $\cP=\varprojlim \cP_l$ is the unipotent
$k$-completion of $G$.  Before proving this, we establish a few facts
about the groups $\cP$ and $\cP_l$.

Recall the filtration $\cP^\bullet$ of $\cP$ defined by
$$\cP^i=\cP\cap (1+\hat{J}_k^i).$$
This induces a filtration on each $\cP_l$:
$$\cP_l=\cP/\cP^{l+1}\supset \cP^2/\cP^{l+1}\supset\cdots\supset
\cP^l/\cP^{l+1}\supset \{1\}.$$
Denote the $i$th term $\cP^i/\cP^{l+1}$ by $L^i$.

Let us determine the structure of the graded algebra
$$\text{Gr}^\bullet \cP = \bop_{i\ge 1} \cP^i/\cP^{i+1}.$$

\begin{lemma}\label{grp} There is an isomorphism
$$\text{\em Gr}^\bullet\cP \stackrel{\cong}{\lra} 
\bop_{i\ge 1} \hat{J}_k^i/\hat{J}_k^{i+1}.$$
\end{lemma}

\begin{proof} For each $i$, define a map $p_i:\cP^i\ra 
\hat{J}_k^i/\hat{J}_k^{i+1}$ by
$$p_i(x)= (x-1) + \hat{J}_k^{i+1}.$$
The map $p_i$ is a homomorphism since
\begin{eqnarray*}
p_i(xy) & = & (xy-1) + \hat{J}_k^{i+1} \\
        & = & (x-1) + (y-1) + (x-1)(y-1) + \hat{J}_k^{i+1} \\
        & = & (x-1) + (y-1) + \hat{J}_k^{i+1} \\
        & = & p_i(x) + p_i(y).
\end{eqnarray*}
The kernel of $p_i$ is precisely the subgroup $\cP^{i+1}$.  We therefore
have an induced injection
$$\overline{p}_i:\cP^i/\cP^{i+1} \lra \hat{J}_k^i/\hat{J}_k^{i+1}.$$
We claim that the map $\overline{p}_i$ is also surjective.  To see this,
consider the composite
$$H_1(G,k)\lra H_1(\cP,k)\lra \cP/\cP^2 \stackrel{\overline{p}_1}{\lra}
\hat{J}_k/\hat{J}_k^2$$
(the map $H_1(\cP,k)\ra \cP/\cP^2$ is the quotient map $(\cP/\Gamma^2\cP)
\otimes k \ra \cP/\cP^2$).  This composite is clearly the canonical 
isomorphism $H_1(G,k)\stackrel{\cong}{\ra}\hat{J}_k/\hat{J}_k^2$.  It
follows that the map $\overline{p}_1$ is surjective.  Now consider the
commutative diagram
$$\xymatrix{
(\cP/\cP^2)^{\otimes i}\ar[r]^{\overline{p}_1^{\otimes i}}\ar[d] &
(\hat{J}_k/\hat{J}_k^2)^{\otimes i} \ar@{>>}[d]^\psi \\
\cP^i/\cP^{i+1}\ar[r]^{\overline{p}_i} & \hat{J}_k^i/\hat{J}_k^{i+1} }$$
Since $\psi$ is surjective and $\overline{p}_i$ is injective, a diagram
chase shows that $\overline{p}_i$ is surjective.
\end{proof}

\begin{lemma}\label{h1gh1pl} The homomorphism $G\ra \cP_l$ induces an
isomorphism $H_1(G,k)\stackrel{\cong}{\ra} H_1(\cP_l,k)$.
\end{lemma}

\begin{proof} Recall the filtration $L^\bullet$ of $\cP_l$ defined by
$L^i=\cP^i/\cP^{l+1}$.  By Lemma \ref{grp}, the graded quotients satisfy
$L^i/L^{i+1}\stackrel{\cong}{\ra} \hat{J}_k^i/\hat{J}_k^{i+1}$.  Hence,
the algebra $\text{Gr}^\bullet_L \cP_l$ is generated by $\text{Gr}^1
\cP_l \cong \cP/\cP^2$.

We have the commutative diagram
$$\xymatrix{
(\cP_l/\Gamma^2\cP_l)^{\otimes i}\ar@{>>}[r]\ar@{>>}[d] & (\cP/\cP^2)^{\otimes
i}\ar@{>>}[d] \\
\Gamma^i\cP_l/\Gamma^{i+1}\cP_l \ar[r]^{\varphi} & L^i/L^{i+1} }$$
The top horizontal arrow is surjective since $\Gamma^2\cP_l \subseteq L^2$,
and the vertical arrows are surjective since the algebras 
$\text{Gr}^\bullet_\Gamma \cP_l$ and $\text{Gr}^\bullet_L\cP_l$ are generated
by $\text{Gr}^1$.  It follows that the map $\varphi$ is surjective.

We have an exact sequence
$$\Gamma^i\cP_l/\Gamma^{i+1}\cP_l \stackrel{\varphi}{\lra} L^i/L^{i+1}
\lra L^i/L^{i+1}\Gamma^i\cP_l \lra 0.$$
Since $\varphi$ is surjective, we see that $L^i=L^{i+1}\Gamma^i\cP_l$.
This implies that $L^i=L^{i+m}\Gamma^i\cP_l$ for all $m\ge 1$.  But since
$L^{l+1}=\{1\}$, we see that $L^i=\Gamma^i\cP_l$ for all $i$.  In particular,
$H_1(\cP_l,k)\cong \hat{J}_k/\hat{J}_k^2$.  Thus we have the commutative
diagram
$$\xymatrix{
H_1(G,k)\ar[r]\ar[dr]_{\cong} & H_1(\cP_l,k)\ar[d]^{\cong} \\
                      & \hat{J}_k/\hat{J}_k^2 }$$
which shows that the map $H_1(G,k)\ra H_1(\cP_l,k)$ is an isomorphism.
\end{proof}

\begin{lemma}  If $H_1(G,k)$ is finite dimensional, then the image of
the map $G\ra\cP_l$ is Zariski dense.
\end{lemma}

\begin{proof}  Denote by ${\mathcal Z}_l$ the Zariski closure of the image
of $G$ in $\cP_l$.  Since the composite
$$H_1(G,k)\lra H_1({\mathcal Z}_l,k)\lra H_1(\cP_l,k)$$
is an isomorphism (Lemma \ref{h1gh1pl}),
the map $H_1({\mathcal Z}_l,k)\ra H_1(\cP_l,k)$ is surjective.  Note that
if $U$ is a unipotent group over $k$, then $H_1(U,\zz) = U/[U,U]$ is
an abelian $k$-group and hence $H_1(U,\zz)\cong H_1(U,k)$.  Lemma
\ref{nilpth1} implies that the map ${\mathcal Z}_l\ra\cP_l$ is surjective;
that is, ${\mathcal Z}_l=\cP_l$. 
\end{proof}

\begin{prop}\label{malcev}  If $H_1(G,k)$ is finite dimensional, then the
map $G\ra\cP$ is the unipotent $k$-completion.
\end{prop}

\begin{proof}  Note that $\cP=\varprojlim \cP_l$ is a prounipotent group over 
$k$ and that the image of $G$ in $\cP$ is Zariski dense.  Denote by $\U$
the unipotent $k$-completion of $G$.  By the universal mapping property,
there is a unique map $\varphi:\U\ra\cP$ making the diagram 
$$\xymatrix{
G\ar[r]\ar[dr] & \U\ar[d]^{\varphi} \\
                      & \cP } $$
commute.  We show that $\varphi$ is an isomorphism.

Since the map $G\ra \cP$ has Zariski dense image, so does the map 
$\U\ra\cP$.  It follows that $\varphi$ is surjective since $\varphi(\U)$
is a closed subgroup of $\cP$ \cite{borel}, p.~47.

To see that $\varphi$ is injective, suppose that $U$ is a unipotent group
over $k$.  Then $U$ is a subgroup of some $U_n\subset GL_n(k)$ for some
$n$ (here, $U_n$ is the subgroup of upper triangular unipotent matrices).
A representation $p:G\ra U$ induces a ring homomorphism $\tilde{p}:
kG\ra {\mathfrak gl}_n(k)$ defined by $\tilde{p}(\sum \alpha_g g) =
\sum \alpha_g p(g)$.  The image of $J_k$ under $\tilde{p}$ lies in
the subalgebra ${\mathfrak n}$ of nilpotent upper triangular matrices.
Thus, the kernel of $\tilde{p}$ contains $J^n_k$ and hence $\tilde{p}$
induces a map
$$\overline{p}:kG/J^n_k\lra {\mathfrak gl}_n(k).$$
Since $\tilde{p}(J_k)\subseteq {\mathfrak n}$, we see that $\cP_{n-1}
\subseteq 1 + J_k/J^n_k$ is contained in $U_n$.  If the image of $G$ is
Zariski dense in $U$, then $\cP_{n-1}\subseteq U$.  Thus, the diagram
$$\xymatrix{
G\ar[r]\ar[dr]_p & \cP_{n-1}\ar[d] \\
                          & U }$$ commutes.  Now, $\U=\varprojlim U_\alpha$,
with each $U_\alpha$ unipotent over $k$.  Applying the above construction
to the compositions $G\ra\U\ra U_\alpha$, we see that $\varphi$ is
injective as follows.  If $\varphi(u)=1$, then $\varphi(u)$ maps to $1$ in
each $\cP_{n-1}$.  But the diagram
$$\xymatrix{
\U\ar[r]\ar[d]_\varphi & U_\alpha  &   &  u\ar[r]\ar[d] & u_\alpha \\
\cP\ar[r]              & \cP_{n-1}\ar[u]   &   &  \varphi(u)=1\ar[r] & 1\ar[u] }$$
shows that $u_\alpha=1$ for each $\alpha$; i.e., $u=1$.
\end{proof}

\begin{cor}  If $k\subset F$ is a field extension, then the map $\U_F
\ra \U_k(F)$ is an isomorphism.
\end{cor}

\begin{proof}  The group $\U_F$ is the set of grouplike elements in 
$FG^\wedge$.  On the other hand, the group $\U_k(F)$ is the set of 
grouplike elements of $kG^\wedge\otimes F=FG^\wedge$.
\end{proof}

\begin{prop}\label{h1uni}  The map $G\ra\cP$ induces an isomorphism
$$H_1(G,k)\stackrel{\cong}{\ra} H_1(\cP,k).$$
\end{prop}

\begin{proof}  Let $A$ be a $k$-vector space with basis $\{e_i\}_{i=1}^n$. 
Then we have a bijection
$$\text{Hom}_{\text{groups}}(G,A) = \text{Hom}_{\text{$k$-vect.~sp.}}
(H_1(G,k),A).$$
On the other hand, if $\sum \alpha_i e_i$ is a vector in $A$, the
map
$$ \sum \alpha_i e_i\mapsto \left(\begin{array}{ccccc}
                     1 & \alpha_1 & \alpha_2 & \cdots & \alpha_n \\ 
                     0 & 1 & 0 & \cdots & 0 \\
                     \vdots & & \ddots & & \vdots \\
                       &  & & 1 & 0 \\
                     0 & 0 & \cdots & 0 & 1 
                 \end{array}\right)$$
identifies $A$ as a unipotent
group over $k$.  The universal mapping property of unipotent 
completion then provides a bijection
$$\text{Hom}_{\text{groups}}(G,A) = \text{Hom}_{\text{$k$-vect.~sp.}}
(H_1(\cP,k),A);$$
that is, $H_1(G,k)$ and $H_1(\cP,k)$ both represent maps of $G$ into
$A$.  It follows that the natural map $H_1(G,k)\ra H_1(\cP,k)$ is an
isomorphism.
\end{proof}

\section{Examples of Unipotent Completions}\label{examples}

\subsection{$G=\zz$}\label{integers}  Let $k$ be a field.  The group
algebra $kG$ is the ring of Laurent polynomials $k[t,t^{-1}]$.  The
augmentation ideal is principal and is generated by $t-1$.  The completion
$kG^\wedge$ is the power series ring $k[[T]]$ with augmentation
$T\mapsto 0$; the ideal $\hat{J}_k$ is the principal ideal $(T)$.  The
canonical map $kG\ra kG^\wedge$ sends $t$ to $1+T$.  If $k$ has
characteristic zero, then $\cP=\{(1+T)^\alpha : \alpha\in k\}$, where
$(1+T)^\alpha = \exp(\alpha\log(1+T))$.  Also, $\cP\cap (1+\hat{J}^l_k)
=\{1\}$ for $l\ge 2$ and hence $\cP=\cP_l$ for all $l$.  Thus, $\cP\cong
k \cong k\otimes_\zz G$, as one would expect.

Let us realize each $\cP_l$ as a group of matrices over $k$.  The
group $\cP_l$ is a subgroup of $\text{Aut}(kG^\wedge/\hat{J}_k^{l+1})$
acting via left multiplication.  In this case $\cP\cong\cP_l$ for
each $l$.  The element $(1+T)^\alpha\in\cP$ is the power series
$$1 + \alpha T + (\frac{\alpha^2}{2}-\frac{\alpha}{2})T^2 +
(\frac{\alpha}{3} -\frac{\alpha^2}{2}+\frac{\alpha^3}{6})T^3 + \cdots .$$
With respect to the basis $1,T,T^2,\dots ,T^l$ of $k[[T]]/(T^{l+1})$, this
element acts via the lower triangular matrix
$$\left(\begin{array}{cccccc}
1 & 0 & 0 &  & \cdots & 0 \\
\alpha & 1 & 0 & & \cdots & 0 \\
c_2(\alpha) & \alpha & 1 & & \cdots & 0 \\
\vdots &  &  &  \ddots & & \vdots \\
c_{l-1}(\alpha) & & & \alpha & 1 & 0 \\
c_l(\alpha) & c_{l-1}(\alpha) & \cdots & c_2(\alpha)  & \alpha & 1
\end{array}\right)$$
where $c_i(\alpha)$ is the coefficient of $T^i$ in the power series 
$(1+T)^\alpha$.  Thus, we have realized $\cP\cong k$ as a closed subgroup
of $\text{Aut}(kG^\wedge/\hat{J}_k^{l+1})$.

If $k=\fp$, then the situation is much different.  Since $(1+T)^p = 1+T^p
\in 1+\hat{J}^p_k$, we see that $1+T$ has finite order in each
$\cP_l = \cP/\cP\cap (1+\hat{J}^{l+1}_k)$.  In fact, $1+T$ has order $p^d$
in $\cP_{p^d},\cP_{p^d+1}, \dots , \cP_{p^{d+1}-1}$.  Note also that
the maps $\zz\ra \cP_l$ are surjective since there are no proper Zariski
dense subgroups of an $\fp$-group (the Zariski topology is discrete).
Each of the homomorphisms
$$\cP_{p^d}\longleftarrow\cP_{p^d+1}\longleftarrow\cdots \longleftarrow
\cP_{p^{d+1}-1}$$ is the identity $\zz/p^d\ra\zz/p^d$.  It follows that
$\cP=\varprojlim \cP_l \cong \zzp$.  Thus, $\cP\cong\zzp\otimes_\zz G$.

Note that this must be correct.  If one hoped that the answer were
$\fp\cong\fp\otimes_\zz G$, then the universal mapping property would 
imply the existence of a nontrivial homomorphism $\fp\ra\zzp$.

Let us make the case $p=2$ explicit.  The group $\cP_l$ is a subgroup of
$\text{Aut}(kG^\wedge/\hat{J}^{l+1}_k)$, acting as left multiplication by
elements of $\cP$.  We shall write down the coordinate rings of the first
few $\cP_l$.  Note that with respect to the basis $1,T,T^2,\dots , T^l$
of $k[[T]]/(T^{l+1})$, the element $1+T$ acts via the lower triangular
$(l+1)\times (l+1)$ matrix
$$\left(\begin{array}{cccccc}
1 & 0 & 0 & \cdots &  & 0  \\
1 & 1 & 0 & \cdots & &  0  \\
0 & 1 & 1 & \cdots &  & 0  \\
\ddots &  &  & \ddots & & \\
      &  &   &   & 1 & 0 \\
0  &  0  &  \cdots &  & 1 & 1 \end{array}\right).$$
For $l=1$, this matrix has order $2$.  The coordinate ring of $\cP_1$ is
$${\mathcal O}(\cP_1) = {\mathbb F}_2[T_{11},T_{12},T_{21},T_{22}]/
(T_{12}, T_{11}-1,T_{22}-1).$$
For $l=2,3$, the above matrix has order $4$.  The coordinate ring of $\cP_2$ is
$${\mathcal O}(\cP_2) = {\mathbb F}_2[T_{ij},1\le i,j \le 3]/
(T_{ij},i<j;T_{ii}-1;T_{21}-T_{32}),$$
and that of $\cP_3$ is
$${\mathcal O}(\cP_3)={\mathbb F}_2[T_{ij},1\le i,j\le 4]/I,$$
where $I$ is the ideal generated by
$T_{ij},i<j;T_{ii}-1;T_{21}-T_{32};T_{32}-T_{43};T_{31}-T_{42};
T_{21}-T_{41}$.
Note that each of these is a two-dimensional ${\mathbb F}_2$-algebra.
For $l=4$, the matrix representing $1+T$ has order $8$.  The coordinate
ring of $\cP_4$ is
$${\mathcal O}(\cP_4) = {\mathbb F}_2[T_{ij},1\le i,j\le 5]/I$$
where $I$ is the ideal generated by
$T_{ij},i<j;T_{ii}-1;T_{21}-T_{41};T_{21}-T_{32};T_{32}-T_{43};
T_{43}-T_{54};T_{31}-T_{42};T_{42}-T_{53};T_{41}-T_{52}$.
It is a three-dimensional ${\mathbb F}_2$-algebra.  To obtain the
${\mathbb F}_{2^d}$-unipotent completion of $\zz$, we take
the inverse limit of the groups $\text{Hom}_{{\mathbb F}_2}
({\mathcal O}(\cP_l),{\mathbb F}_{2^d})$.

\subsection{$G=\zz/n$}\label{zmodn}  Here the group algebra is
$kG=k[t]/(t^n-1)$; the augmentation ideal is $J_k=(t-1)$.  If $k$ has
characteristic zero, then since
\begin{eqnarray*}
(t-1)^2 & = & t^2-2t+1 \\
            & = & t^n+t^2-2t \\
            & = & t(t^{n-2}+\cdots +t+2)(t-1),
\end{eqnarray*}
and since both $t$ and $t^{n-2}+\cdots +t+2$ are units in $kG$, we have
$J^2_k=J_k$.  Thus, the completion $kG^\wedge$ is the algebra 
$\varprojlim kG/J_k^l = k$.  The only grouplike element is $1$, and hence
the unipotent $k$-completion is trivial.

If $k=\fp$, then we must consider separately two cases.  First assume that $p$ is prime to $n$.  In this case, the above factorization shows
again that $kG^\wedge=k$ and that $\cP=1$.  If $p$ divides $n$, then we
may assume that $n=p^d$ for some $d\ge 1$ (note that the completion of
$G\times H$ is the product of the completions of $G$ and $H$).  In this
case we have $J^{p^d}=0$ since $(1-t)^{p^d} = (1-t^{p^d})=(1-t^n)=0$.
Thus, $kG^\wedge \cong kG$ and $\cP$ is the collection of the $t^i$, $0\le i\le
p^d$.  This group is $\zz/p^d$.

In the characteristic zero case, then, we see that the completion of
$\zz/n$ is $\zz/n\otimes_\zz k$ and in the characteristic $p$ case it
is $\zz/n\otimes_\zz \zzp$.  So, if $G$ is a finitely generated abelian
group, the $\zq$-completion of $G$ is $G\otimes_\zz\zq$ and the
$\fp$-completion is $G\otimes_\zz\zzp$.

\subsection{$p$-completion v.~$\fp$-completion}\label{pisfp}  Suppose
that $H_1(G,\fp)$ is finite dimensional.  Then the $\fp$-completion of
$G$ is the $p$-completion.  Indeed, the $p$-completion of $G$ is the 
inverse limit
$$G^{\wedge p} = \varprojlim G/\Gamma^r_pG,$$ where $\Gamma^\bullet_p G$ is the $p$-lower central series of $G$.  Each
group $G/\Gamma^r_pG$ is a finite $p$-group if $H_1(G,\fp)$ is
finite dimensional.  On the other hand, each group $\cP_l$ is a unipotent
group over $\fp$.  Such groups are nilpotent and since $\cP_l$ is a closed
subgroup of some $U_n(\fp)$, it is finite.  Thus, we have a homomorphism
$$G^{\wedge p}\lra \cP.$$  But since $G^{\wedge p}$ is a prounipotent
group over $\fp$, we have a unique map
$$\cP\lra G^{\wedge p},$$ and it is clear that these maps are inverse
isomorphisms.  Alternatively, both $G^{\wedge p}$ and $\cP$ satisfy the
same universal mapping property and are therefore isomorphic via a
unique isomorphism.

\subsection{Coordinate rings of curves}\label{coordrings}  Let $X$ be
a smooth affine curve over $k$ and let $x$ be a closed point on $X$.  
Denote by $\Ox$ the local ring at $x$; it has maximal ideal $\mx$.  
We have a split extension
$$1\lra K(\Ox)\lra \slno \stackrel{\mod\mx}{\lra}\slnk\lra 1.$$
The group $K(\Ox)$ is filtered by powers of the maximal ideal
$\mx$; denote this filtration by $K^\bullet(\Ox)$.  By a theorem of
W.~Klingenberg \cite{kling}, $K^\bullet(\Ox)$ is the lower central series
of $K(\Ox)$ provided $\text{char}(k)\ne 2$ or $k\ne {\mathbb F}_3$.
The graded quotients satisfy
$$K^l/K^{l+1}\cong \lslnk;$$
in particular, $H_1(K(\Ox),\zz)\cong\lslnk$.

This is a special case of the following.

\begin{prop}\label{centralk}  Let $G$ be a group such that the graded
quotients $\Gamma^lG/\Gamma^{l+1}G$ are finite dimensional 
$k$-vector spaces.  Let $\hat{G}=\varprojlim G/\Gamma^lG$.  
Then the map $G\ra\hat{G}$ is the unipotent $k$-completion.
\end{prop}

\begin{proof} This is essentially the same as Proposition \ref{malcev}.
Each $G/\Gamma^lG$ is a unipotent $k$-group by hypothesis and since $G$
surjects onto each $G/\Gamma^lG$, the image of $G$ in $\hat{G}$ is
Zariski dense.  Since $H_1(G,k)$ is finite dimensional, the unipotent
$k$-completion of $G$ is the group $\cP\subset kG^\wedge$.  The universal
mapping property provides a unique map $\Phi:\cP\ra\hat{G}$ of prounipotent
$k$-groups; it is surjective since $\Phi(\cP)$ is a closed Zariski dense
subgroup of $\hat{G}$.  The injectivity of $\Phi$ follows since the
composition $G\ra\cP\ra\cP_l$ must factor through some $G/\Gamma^lG$.
\end{proof}

Now, for the group $K(\Ox)$, $\text{char}(k)\ne 2, k\ne {\mathbb F}_3$,
the lower central series satisfies the hypotheses of the proposition. Thus,
the unipotent $k$-completion is the group
$$K(\Ohat) = \text{ker}\{\slnohat\lra\slnk\},$$ where $\Ohat$ is the
$\mx$-adic completion of $\Ox$.

\subsection{$G=\slnz$}\label{slnz}  We have the exact sequence
$$1\lra \Gamma(n,p)\lra\slnz\lra\slnfp\lra 1.$$  If $n\ge 3$, the 
filtration of $\Gamma(n,p)$ by powers of $p$ is the lower central series
\cite{bms}.  Since the graded quotients 
$\Gamma(n,p^l)/\Gamma(n,p^{l+1})$ are isomorphic to ${\mathfrak sl}_n
(\fp)$, Proposition \ref{centralk} implies that the $\fp$-completion of 
$\Gamma(n,p)$ is the group
$$\hat{\Gamma} = \varprojlim \Gamma(n,p)/\Gamma(n,p^l).$$  This group
fits into the exact sequence
$$1\lra\hat{\Gamma}\lra SL_n(\zzp)\lra\slnfp\lra 1.$$
Note that the $\fp$-completion of $\slnz$ is trivial.  Also, if $\ell\ne p$
the $\fl$-completion of $\Gamma(n,p)$ is trivial since $H_1(\Gamma(n,p),\fl)
=0$.

If $k$ is any field of characteristic $p$, then the unipotent $k$-completion of $\Gamma(n,p)$ is obtained by taking the $k$-form of
$\hat{\Gamma}$.  This is constructed as follows.  The first group in the
inverse system is $U_2=\lslnk$.  The next is a nonsplit extension
$$0\lra\lslnk\lra U_3\lra\lslnk\lra 0.$$  The group $U_3$ is the unipotent
group defined by the diagram
$$\xymatrix{
0 \ar[r] & {\mathfrak sl}_n(\fp)\ar[d]\ar[r] & \Gamma(n,p)/\Gamma(n,p^3)
\ar[d]\ar[r] & {\mathfrak sl}_n(\fp)\ar[d]\ar[r] & 0 \\
0\ar[r] & \lslnk\ar[r] & U_3\ar[r] & \lslnk\ar[r] & 0. }$$
In general, we have a diagram of central extensions
$$\xymatrix{
0 \ar[r] & {\mathfrak sl}_n(\fp)\ar[d]\ar[r] & \Gamma(n,p)/\Gamma(n,p^r)
\ar[d]\ar[r] & \Gamma(n,p)/\Gamma(n,p^{r-1})\ar[d]\ar[r] & 1 \\
0\ar[r] & \lslnk\ar[r] & U_r\ar[r] & U_{r-1}\ar[r] & 1 }$$
and the $k$-completion of $\Gamma(n,p)$ is $\U_k=\varprojlim U_r$.

\section{Relative Completions}\label{relative} 
Let $S$ be a reductive algebraic group over $k$ and let $\rho:\Gamma
\ra S$ be a representation with Zariski dense image. Denote by $\G$
the completion of $\Gamma$ with respect to $\rho$; it is an extension
$$1\lra\U\lra\G\lra S\lra 1.$$  In this section we construct some
examples of relative completions.  We first prove a few basic facts.

First, let $A$ be a rational representation of $S$ and suppose that we
have an extension in the category of algebraic groups over $k$:
$$0\lra A\lra G\lra S\lra 1.$$  Then this extension splits \cite{borel},
p.~158.  This extends to unipotent kernels; a proof of the following may
be found in \cite{hain1}, Prop. 4.3.

\begin{prop}  If $S$ is a reductive algebraic group over $k$ and if
$$1\lra U\lra G\lra S\lra 1$$ is an extension in the category of algebraic
groups over $k$, where $U$ is unipotent, then the extension splits.
\end{prop}

From this, one easily deduces that the relative completion
$$1\lra\U\lra\G\lra S\lra 1$$ is a split extension.

Denote the kernel of $\rho$ by $T$ and the unipotent $k$-completion of
$T$ by $\T$.  Then the map $T\ra\U$ induces a unique map $\Phi:\T\ra\U$.
Denote by $\K$ the kernel of $\Phi$ and by $L$ the image of $\rho$.

\begin{prop} Suppose that $H_1(T,k)$ is finite dimensional.  If the action
of $L$ on $H_1(T,k)$ extends to a rational representation of $S$
(e.g., $L=S$), then $\K$ is central in $\T$.
\end{prop}

\begin{proof} \cite{hain1}  The group $\Gamma$ acts on the completion 
$kT^\wedge$ by conjugation.  This action preserves the filtration by
powers of $\hat{J}$, so $\Gamma$ acts on the associated graded algebra
$$\text{Gr}^\bullet_{\hat{J}} kT^\wedge = \bop_{m\ge 0} \hat{J}^m/
\hat{J}^{m+1}.$$
If $H_1(T,k)$ is finite dimensional, then each algebra $kT^\wedge/
\hat{J}^l$ is finite dimensional.  Thus, each of the groups $\text{Aut}
(kT/J^l)$ is an algebraic group.  Since $\text{Gr}^\bullet_{\hat{J}} 
kT^\wedge$ is generated by $J/J^2=H_1(T,k)$, it follows that
$\text{Aut}(kT^\wedge)$, the group of augmentation preserving algebra
automorphisms of $kT^\wedge$, is a proalgebraic group which is an 
extension of a subgroup of $\text{Aut}H_1(T,k)$ by a prounipotent group
${\mathcal V}$:
\begin{equation}\label{extension}
1\lra {\mathcal V}\lra \text{Aut}(kT^\wedge)\lra \text{Aut}H_1(T,k).
\end{equation}

If the action of $\Gamma$ on $H_1(T,k)$ factors through a rational
representation $S\ra\text{Aut}H_1(T,k)$, then we can form a proalgebraic
group extension
$$1\lra {\mathcal V}\lra {\mathcal E}\lra S\lra 1$$ of
$S$ by the prounipotent group ${\mathcal V}$ by pulling back
the extension (\ref{extension}) along $S\ra \text{Aut}H_1(T,k)$. Since
the map $\Gamma\ra\text{Aut}H_1(T,k)$ factors through $S$, we can lift the
map $\Gamma\ra\text{Aut}(kT^\wedge)$ to a map $\Gamma\ra {\mathcal E}$ whose
composition with ${\mathcal E}\ra S$ is $\rho:\Gamma\ra S$.
By the universal mapping property of relative completion, this induces a 
homomorphism $\G\ra {\mathcal E}$. Since the composite
$$\T\lra\U\lra\G\lra {\mathcal E}\lra \text{Aut}(kT^\wedge)$$ is the action
of $\T\subset kT^\wedge$ on $kT^\wedge$ by inner automorphisms, we see
that the kernel of this map is the center of $\T$.  It follows that the
kernel of $\Phi$ is central in $\T$.
\end{proof}

In \cite{hain1}, Hain states the result (due to Deligne) that if $H^1(L,A)
=0$ for all rational representations $A$ of $S$, then $\Phi:\T\ra\U$ is
surjective.  This works well in characteristic zero thanks to the vanishing
results for arithmetic groups due to Ragunathan \cite{ragu}.  In positive
characteristic, however, this criterion is not very useful.  We present here
a simplified version which suits our purposes.

\begin{prop}\label{phisurj} Assume that $H_1(T,k)$ is finite dimensional.
If $\rho:\Gamma\ra S$ is surjective, then $\Phi:\T\ra\U$ is surjective.
\end{prop}

\begin{proof}  We first show that the composition
$$H_1(T,k)\lra H_1(\T,k)\stackrel{\Phi_*}{\lra} H_1(\U,k)$$ is surjective.  
Let $A$ be the cokernel of this
map; it is a rational representation of $S$.  Pushing out the extension
$$1\lra\U\lra\G\lra S\lra 1$$ along the map $\U\ra H_1(\U,k)\ra A$,
we obtain an extension of algebraic groups
$$0\lra A\lra G\lra S\lra 1.$$  
This splits; say $s:S\ra G$ is a splitting.  Since
the image of $\Gamma$ in $\G$ is Zariski dense, the image of $\Gamma$
in $G$ is Zariski dense.  But since the diagram
$$\xymatrix{
\Gamma\ar[r]\ar[dr]_\rho & G \\
       & S\ar[u]_s }$$
commutes, we see that the image of $s$ is Zariski dense in $G$.  This
forces $A$ to be trivial.

We have shown that the map $\Phi_*:H_1(\T,k)\ra H_1(\U,k)$ is surjective.
In particular, both of these vector spaces are finite dimensional
($H_1(\T,k)$ is finite dimensional by Proposition \ref{h1uni}).  Also,
both $\T$ and $\U$ are their own unipotent $k$-completions.  Thus we have
isomorphisms \linebreak[4] $\T\stackrel{\cong}{\ra}\cP(k\T^\wedge)$ and
$\U\stackrel{\cong}{\ra}\cP(k\U^\wedge)$, where $\cP(k{\mathcal H}^\wedge)$ is
the set of grouplike elements of $k{\mathcal H}^\wedge$.  Now, $H_1(\T,k)\cong
\hat{J}_{\T}/\hat{J}_{\T}^2$ and $H_1(\U,k)\cong \hat{J}_{\U}/\hat{J}_{\U}^2$,
where $\hat{J}_{\mathcal H}$ is the augmentation ideal in 
$k{\mathcal H}^\wedge$.  Since $\Phi_*$ is surjective, we see that the
map $$\text{Gr}^\bullet\hat{\Phi}:\text{Gr}^\bullet k\T^\wedge\lra
\text{Gr}^\bullet k\U^\wedge$$ is surjective.  But then the induced
map $\hat{\Phi}:k\T^\wedge\ra k\U^\wedge$ is surjective (\cite{quillen},
p.~266).  This implies that $\Phi:\T\ra\U$ is surjective since
$\hat{\Phi}$ maps grouplike elements to grouplike elements.
\end{proof}

We are now able to compute some examples of relative completions.

\subsection{}\label{localrel}  Let $X$ be a smooth affine curve over $k$
and let $\Ox$ be the local ring at a closed point $x$ with maximal ideal
$\mx$.  In the previous section we showed that the unipotent 
$k$-completion of $K(\Ox)=\{A\in \slno : A\equiv I\mod\mx\}$ is the
group $K(\Ohat)$, provided that $\text{char}(k)\ne 2$ or $k\ne {\mathbb
F}_3$.  Let $\rho:\slno\ra\slnk$ be reduction modulo $\mx$ and consider
the completion relative to $\rho$
$$1\lra \U\lra\G\lra\slnk\lra 1.$$
There is a surjective map $\Phi:K(\Ohat)\ra\U$ with central kernel.  An
easy calculation shows that the center of $K(\Ohat)$ is trivial.  Since
the above extension splits, we see that the group $\G$ is the group
$\slnohat$.

\subsection{}\label{ktrel}  Consider the curve $X={\mathbb A}_k^1$.  The
coordinate ring is the polynomial algebra $k[t]$.  If $k$ is a finite field
or a number field, then the filtration of 
$$K(k[t])=\text{ker}\{\slnkt\stackrel{\rho}{\lra}\slnk\}$$ by powers of
$t$ is the lower central series, provided that $n\ge 3$ \cite{knudson3}.
Thus, the unipotent $k$-completion of $K(k[t])$ is the group $K(k[[T]])$.
Since $\rho:\slnkt\ra\slnk$ is surjective, and since $K(k[[T]])$ has
trivial center, we see that the relative completion is $\slnpow$.

\subsection{}\label{truncrel}  Let $k$ be any field and consider the group
$SL_n(k[t]/t^l)$.  Reduction modulo $t$ gives a representation
$\rho:SL_n(k[t]/t^l)\ra\slnk$.  The kernel $K$ is filtered by powers of $t$;
this filtration is finite in length and is the lower central series.  The
graded quotients are all isomorphic to $\lslnk$.  It follows that the
unipotent $k$-completion of $K$ is
$$\varprojlim K/K^l \cong K$$ (this holds for all $n$).  Thus, the relative
completion is simply the group $SL_n(k[t]/t^l)$.

\subsection{}\label{zrel}  Consider the extension
$$1\lra\Gamma(n,p)\lra\slnz\stackrel{\rho}{\lra}\slnfp\lra 1$$
for $n\ge 3$.  The $\fp$-completion of $\Gamma(n,p)$ is the group
$\hat{\Gamma}$:
$$1\lra\hat{\Gamma}\lra SL_n(\zzp)\lra\slnfp\lra 1.$$
Denote the completion of $\rho:\slnz\ra\slnfp$ by $\G$ and let $\U$ be
its prounipotent radical.  The map $\Phi:\hat{\Gamma}\ra\U$ is surjective
since $\rho$ is and its kernel is central in $\hat{\Gamma}$.  Clearly,
though, the center of $\hat{\Gamma}$ is trivial and hence $\Phi$ is
an isomorphism.  Thus, $\G$ is the semidirect product
$$\G\cong \hat{\Gamma}\cdot\slnfp.$$  (Recall that the relative
completion is a split extension.)  This group is {\em not} isomorphic
to $SL_n(\zzp)$ since the projection map $SL_n(\zzp)\ra\slnfp$
does not split.  Note also that $SL_n(\zzp)$ is not proalgebraic over $\fp$
since the rings $\zz/p^d$ are not $\fp$-algebras.

Thus, the group $\G$ is a bit mysterious in the sense that we do not have
an explicit description of it in terms of matrices, nor do we have a 
formula for the homomorphism $\slnz\ra\G$.

\subsection{}\label{galois}  Let $\Gamma$ be the absolute Galois group
of the finite field $\fp$.  The group $\Gamma$ is isomorphic to $\hat{\zz}$, the profinite completion of $\zz$.  Let $\ell$ be a prime
different from $p$ and consider the action of $\Gamma$ on the $\ell$th
roots of unity in $\fpbar$.  This defines a homomorphism $\rho:\Gamma
\ra\fl^\times$ as follows.  Let $\zeta_\ell$ be a primitive $\ell$th root
of unity and let $\sigma$ be an element of $\Gamma$.  Write
$$\sigma(\zeta_\ell) = \zeta_\ell^{\rho(\sigma)}, \qquad \rho(\sigma)\in
\{1,2,\dots ,\ell-1\}.$$  (Note that $\rho(\sigma)\ne 0$ since $\sigma(1)
=1$.)  This defines the map $\rho$.  Since $\rho$ maps the Frobenius
automorphism to a generator of $\fl^\times$, we see that $\rho$ is
surjective.  The kernel of $\rho$ consists of those automorphisms that
fix $\fp(\zeta_\ell)$; that is, $\text{ker} \rho = \text{Gal}(\fpbar/
\fp(\zeta_\ell))$.  This group is also isomorphic to $\hat{\zz}$ and the
inclusion $\text{ker} \rho\ra\Gamma$ is simply multiplication by 
$\ell-1$.  Note that $\hat{\zz}\cong \prod_q \zz_q$.  Thus, we have an
extension
$$0\lra\hat{\zz}\stackrel{\ell-1}{\lra}\hat{\zz}\stackrel{\rho}{\lra}
\fl^\times\lra 1.$$  The $\fl$-completion of $\hat{\zz}$ is clearly
$\zz_\ell$ and since $\rho$ is surjective we have a surjection $\Phi:
\zz_\ell\ra\U$, where $\U$ is the prounipotent radical of the completion
of $\Gamma$ relative to $\rho$.  The kernel of $\Phi$ is central in 
$\zz_\ell$, but that is not a useful piece of information in this case.
We claim that $\Phi$ is injective so that the relative completion is the
semidirect product
$$\xymatrix{
0\ar[r] &\zz_\ell\ar[r] &\G\ar[r] &\fl^\times\ar[r]\ar@/^/[l] & 1,}$$ 
where $\fl^\times$ acts via
multiplication on $\zz_\ell$.

Recall that in the general situation we denote by $L$ the image of $\rho:
\Gamma\ra S$.  Assume that $H^1(L,A)$ vanishes for all rational
representations $A$ of $S$ and that $H^2(L,A)$ vanishes for all nontrivial
$A$.  Consider the extension
\begin{equation}\label{centralext}
0\lra\K\lra G\lra L\lra 1
\end{equation}
where $\K$ is the kernel of $\Phi:\T\ra\U$.  Since $H_1(L,k)=0$, there is
a universal central extension with kernel a $k$-vector space
$$0\lra H_2(L,k)\lra \tilde{L}\lra L\lra 1.$$ Its cocycle is the identity
map in 
$$\text{Hom}(H_2(L,k),H_2(L,k))\cong H^2(L,H_2(L,k)).$$  The extension
(\ref{centralext}) is classified by a linear map $\psi:H_2(L,k)\ra\K$.
Under the above stated conditions on $H^1(L,A)$ and $H^2(L,A)$, the map
$\psi$ is surjective \cite{hain1}, Prop. 4.13.

In the example under consideration, we have $L=S=\fl^\times$.  If $A$
is a rational $S$-module, then $H^1(S,A)=0=H^2(S,A)$ since $|S|=\ell-1$
is invertible in $A$ ($A$ is an $\fl$-vector space).  Since $H_2(\fl^\times,
\fl)=0$, we see that $\K=0$.

\section{Cohomology}\label{cohomology}

Suppose that $\pi$ is a projective limit of groups, $\pi=\varprojlim 
\pi_\alpha$, and let $k$ be a field.  We define the continuous cohomology
of $\pi$ to be $$H^i_{\text{cts}}(\pi,k) = \varinjlim H^i(\pi_\alpha,k).$$
For example if $\pi$ is the Galois group of a field extension $L/F$,
then $H^i_{\text{cts}}(\pi,k)$ is simply the usual Galois cohomology.
Note that there is a natural map $H^i_{\text{cts}}(\pi,k)\ra H^i(\pi,k)$.
It is obvious that $H^0_{\text{cts}}(\pi,k)=H^0(\pi,k)$.

\begin{lemma}\label{ctsh2}  Let $G$ be a group and let $\U$ be its
unipotent $k$-completion.  Then the map 
$$\mhcts(\U,k)\lra H^2(G,k)$$ is injective.  If $H_1(G,k)$ is finite
dimensional, then the map $$H^1_{\text{\em cts}}(\U,k)\lra H^1(G,k)$$
is an isomorphism.
\end{lemma}

\begin{proof}  Let $\alpha\in\hcts(\U,k)$.  This class corresponds to a 
central extension
$$0\lra k\lra\U_\alpha\lra\U\lra 1$$ in the category of proalgebraic
$k$-groups.  Observe that $\U_\alpha$ is prounipotent.  Suppose that
$\alpha$ restricts to $0$ in $H^2(G,k)$.  Then we have a commutative 
diagram of extensions
$$\xymatrix{
0\ar[r] & k\ar[r] & \U_\alpha\ar[r] & \U\ar[r] & 1 \\
0\ar[r] & k\ar[u]_{=}\ar[r] & G_\alpha\ar[u]\ar[r]  & G\ar[u]\ar[r] & 1. }$$
Since $\alpha$ maps to $0$ in $H^2(G,k)$, the bottom extension splits.
We then have the composite map $G\ra G_\alpha\ra\U_\alpha$.  By the
universal mapping property of $\U$, we get a unique map $\U\ra\U_\alpha$ making the diagram
$$\xymatrix{
G\ar[r]\ar[dr] & \U\ar[d] \\
                        & \U_\alpha }$$
commute; i.e., the top extension splits.  Thus, $\alpha=0$.

If $H_1(G,k)$ is finite dimensional, then the unipotent $k$-completion
of $G$ is the group $\cP\subset kG^\wedge$.  By Lemma \ref{h1gh1pl},
the map $H_1(G,k)\ra H_1(\cP_l,k)$ is an isomorphism for each $l$.
It follows that the map
$$H^1_{\text{cts}}(\cP,k)=\varinjlim H^1(\cP_l,k)\lra H^1(G,k)$$
is an isomorphism.
\end{proof}

\noindent {\em Remark.}  The above argument does not work to show that
$H^2(\U,k)$ injects into $H^2(G,k)$.  An element $\alpha$ corresponds
to an extension of $\U$ by $k$, but this may not be an extension of
proalgebraic groups.  A class in $\hcts(\U,k)$ actually corresponds to
a compatible sequence of extensions of algebraic groups; this allows
us to use the universal mapping property of $\U$.

\begin{cor}\label{h1iso} The map $H^1_{\text{\em cts}}(\U,k)\ra H^1(\U,k)$
is an isomorphism.
\end{cor}

\begin{proof} We have the commutative diagram
$$\xymatrix{
H^1_{\text{cts}}(\U,k)\ar[r]\ar[dr]_{\cong} & H^1(\U,k)\ar[d] \\
        & H^1(G,k).}$$
The vertical arrow is an isomorphism by Proposition \ref{h1uni}.
\end{proof}

\begin{cor}\label{h2inj}  The map $\mhcts(\U,k)\ra H^2(\U,k)$ is
injective. \hfill $\qed$
\end{cor}

\begin{cor}\label{relinj}  Let $\rho:\Gamma\ra S$ be a split surjective
representation and let $\G$ be the completion relative to $\rho$.  Assume
that $H_1(T,k)$ is finite dimensional {\em(}$T=\text{\em ker} \rho${\em)}
and that $\Phi:\T\ra\U$ is an isomorphism.  Then the restriction map
$$\mhcts(\G,k)\lra H^2(\Gamma,k)$$ is injective.
\end{cor}

\begin{proof} Note that $\G$ is the inverse limit of the groups $\G_\alpha$
with each $\G_\alpha$ an extension
$$1\lra \U_\alpha\lra \G_\alpha\lra S\lra 1.$$
For each $\alpha$, we have a Hochschild--Serre spectral sequence
$$E_2^{i,j}(\G_\alpha) = H^i(S,H^j(\U_\alpha,k))\Longrightarrow H^{i+j}
(\G_\alpha,k).$$ Taking the direct limit of these sequences, we get a 
spectral sequence $E(\G)$ satisfying
$$E_2^{i,j}(\G) = H^i(S,\varinjlim H^j(\U_\alpha,k)) \Longrightarrow
\varinjlim H^{i+j}(\G_\alpha,k);$$
that is,
$$E_2^{i,j}(\G) = H^i(S,H^j_{\text{cts}}(\U,k))\Longrightarrow
H^{i+j}_{\text{cts}}(\G,k).$$
Since $\Phi:\T\ra\U$ is an isomorphism we have $H^1_{\text{cts}}(\U,k)
\cong H^1(T,k)$ (Lemma \ref{ctsh2}). 
Since the extension $$1\lra T\lra\Gamma\lra S\lra 1$$ splits, we have
$$H^i(S,H^j_{\text{cts}}(\U,k)) \stackrel{\cong}{\lra} H^i(S,H^j(T,k))$$
for $j=0,1$; that is, we have an isomorphism
$$E_2^{i,j}(\G)\stackrel{\cong}{\lra}E_2^{i,j}(\Gamma), \qquad j=0,1,$$ 
where $E(\Gamma)$ is
the Hochschild--Serre spectral sequence for the above extension. Moreover, the 
differential $d^2:E_2^{i,1}\ra E_2^{i+2,0}$ vanishes for both spectral
sequences since the extensions are split.
By Lemma \ref{ctsh2}, there is an inclusion of $S$-modules
$\hcts(\U,k)\ra H^2(T,k)$ and hence $E_2^{0,2}(\G)$ injects into
$E_2^{0,2}(\Gamma)$.  These facts together imply the following:
\begin{enumerate}
\item $E_\infty^{2,0}(\G)=E_\infty^{2,0}(\Gamma)$;
\item $E_\infty^{1,1}(\G) = E_\infty^{1,1}(\Gamma)$;
\item $E_\infty^{0,2}(\G)\hookrightarrow E_\infty^{0,2}(\Gamma)$.
\end{enumerate}
It follows that $\hcts(\G,k)$ injects into $H^2(\Gamma,k)$.
\end{proof}

\subsection{} Let $X$ be a smooth affine curve over $k$ and let $\Ox$
be the local ring at $x\in X$.  In the previous section, we showed that the
completion of $\slno\ra\slnk$ is the extension
$$1\lra K(\Ohat)\lra\slnohat\lra\slnk\lra 1,$$ provided that $\text{char}
(k)\ne 2$ or $k\ne {\mathbb F}_3$ (recall that we need these restrictions on
$k$ to guarantee that $K(\Ohat)$ is the unipotent $k$-completion of
$K(\Ox)$).  By Corollary \ref{relinj}, we see that the map
$$\hcts(\slnohat,k)\lra H^2(\slno,k)$$ is injective for all $n$.

\subsection{}  Let $k$ be a finite field or a number field and consider the
completion of $\slnkt\ra\slnk$ for $n\ge 3$.  This is the split extension
$$1\lra K(k[[T]])\lra\slnpow\lra\slnk\lra 1.$$

\begin{theorem}\label{powercoho}  Let $k$ be a finite field or a number
field.  Then for all $n\ge 3$, $\mhcts(\slnpow,k)=0$.
\end{theorem}

\begin{proof}  We have an injection
$$\hcts(\slnpow,k)\lra H^2(\slnkt,k).$$  If $k$ is a number field then
by \cite{knudson1}, the group $H^2(\slnkt,k)$ coincides with $H^2(
\slnk,k)$.  We then have the chain of equalities
\begin{eqnarray*}
H^2(\slnkt,k) & = & H^2(\slnk,k) \\
                     & = & H_2(\slnk,k) \\
                     & = & H_2(\slnk,\zz)\otimes k \\
                     & = & K_2(k)\otimes k.
\end{eqnarray*}
Since $K_2(k)$ is a torsion group when $k$ is a number field, we see that
$H^2(\slnkt,k)=0$.

If $k$ is a finite field, then we use the following chain of equalities
\begin{eqnarray*}
H^2(\slnkt,k) & = & H_2(\slnkt,k) \\
                     & = & H_2(\slnkt,\zz)\otimes k \\
                     & = & K_2(k[t])\otimes k \\
                      & = & K_2(k)\otimes k,
\end{eqnarray*}
where the last equality is the fundamental theorem of algebraic 
$K$-theory.  Since $K_2(k)=0$ for finite fields, we are done.
\end{proof}

\noindent {\em Remark.}  It is easy to see that if $E$ is a field with
$\text{char}(E)\ne\text{char}(k)$, then $H^\bullet_{\text{cts}}(
\slnpow,E)=H^\bullet(\slnk,E)$.  Also, Theorem \ref{powercoho} could be
deduced for $n\ge 5$ using homological stability for $SL_n$ together
with existing calculations of the groups $K_2(k[t]/t^l)$.  The above
proof avoids this (and works for $n\ge 3$).

\medskip

In \cite{hain2}, Hain calls a prounipotent group $\cP$ pseudonilpotent
if $$H^\bullet_{\text{cts}}(\cP,k)\stackrel{\cong}{\lra} H^\bullet(\cP,k)$$ 
(his definition is
actually more general).  Let $k$ be a number field and consider the
group $K(k[[T]])$.  We know that $\hcts(K(k[[T]]),k)$ injects into
$H^2(K(k[[T]]),k)$.   Suppose that this map is an isomorphism.  Then we
would have an isomorphism
$$H^2(\slnpow,k)\stackrel{\cong}{\lra} H^2(\slnkt,k)\stackrel{\cong}
{\lra} H^2(\slnk,k).$$  Dualizing and taking the direct limit over $n$,
we would have an isomorphism
$$K_2(k[[T]])\otimes k \stackrel{\cong}{\lra} K_2(k)\otimes k.$$
Gabber's rigidity theorem \cite{gabber} asserts that $K_2(k[[T]],
\zz/n)\cong K_2(k,\zz/n)$.  We would therefore have an isomorphism
$K_2(k[[T]])\cong K_2(k)$.  This is not the case, however, since the
kernel of the map $K_2(k[[T]])\ra K_2(k)$ is the relative group
$K_2(k[[T]],(T))$ and this group is nontrivial \cite{vdk} (indeed, it
is uniquely divisible).  It follows that the group $K(k[[T]])$ is not
pseudonilpotent.

\subsection{} Consider the commutative diagram of extensions
$$\xymatrix{
1\ar[r] & \Gamma(n,p)\ar[d]\ar[r] & \slnz\ar[d]\ar[r] & \slnfp\ar[d]\ar[r]
& 1 \\
1\ar[r] & \hat{\Gamma}\ar[r] & SL_n(\zzp)\ar[r] & \slnfp\ar[r] & 1.} $$
We showed (\ref{slnz}) that $\hat{\Gamma}$ is the $\fp$-completion
of $\Gamma(n,p)$ if $n\ge 3$.  While Corollary \ref{relinj} does not
apply in this case, we can still use the argument in the proof to deduce
that the map 
$$\hcts(SL_n(\zzp),\fp)\lra H^2(\slnz,\fp)$$ is injective for $n\ge 3$.
Note that $H^2(\slnz,\fp)\cong H_2(\slnz,\fp)$.  The group $H_2(\slnz,\zz)$ is $2$-torsion for all $n\ge 3$; it is equal to $\zz/2
\oplus\zz/2$ for $n=3,4$ \cite{vdk1}, and to $\zz/2$ for $n\ge 5$
\cite{milnor}.  It follows that if $p\ge 3$, then $H_2(\slnz,\fp)=0$.

\begin{theorem} Let $n\ge 3$.  If $p\ge 3$, then $\mhcts(SL_n(\zzp),\fp)=0$.
When $p=2$ we have
$$\dim_{{\mathbb F}_2}\mhcts(SL_n({\mathbb Z}_2),{\mathbb F}_2) \le
\begin{cases}
       2 & n=3,4 \\
       1 & n\ge 5.
\end{cases}$$
\end{theorem}

\end{document}